\documentstyle[amssymb]{amsart}

\newcommand	{\comment}[1]	{}
\newcommand{\mute}[2] {}
\newcommand	{\printname}[1]	{}
\newcommand{\labell}[1]	{\label{#1}\printname{#1}}
\newcommand{\bibitemm}[2]{\bibitem[#1]{#2}\printname{#2}}
\numberwithin{equation}{section}
\newtheorem {Theorem}		{Theorem}
\newtheorem {Lemma}     	{Lemma}[section]

\newtheorem {Corollary} [Lemma]	{Corollary}
\newtheorem {Proposition} [Lemma]{Proposition}

\theoremstyle{definition}
\newtheorem{Definition}[Lemma]{Definition}
\newtheorem{Notation}[Lemma]{Notation}

\theoremstyle{remark}
\newtheorem{Remark}[Lemma]{Remark}

\newtheorem{Example}[Lemma]{Example}

\def	\R	{{\Bbb R}}
\def	\C	{{\Bbb C}}
\def	\Z	{{\Bbb Z}}

\def	\PP	{{\Bbb P}}

\def    \t      {{\frak  t}}
\def    \g      {{\frak  g}}
\def    \h      {{\frak  h}}

\def    \fg      {{\frak  g}}
\def    \fh      {{\frak  h}}

\newcommand 	{\dd}[1]	{{\partial \over \partial #1}}
\def	\ssminus{\smallsetminus}
\def	\to	{\longrightarrow}

\def	\Cinf	{C^\infty}

\def	\ol	{\overline}

\def	\tPhi	{{\tilde{\Phi}}}
\def	\inv	{^{-1}}
\def	\calV	{{\cal V}}
\def	\half	{{1 \over 2}}

\def	\Lie	{{\mbox{Lie}}}
\def	\reg	{{\text{reg}}}
\def	\loc	{{\text{loc}}}

\begin{document}
\title[The centralizer of invariant functions] {The centralizer of
invariant functions and division properties of the moment map}

\author[Y. Karshon]{Yael Karshon}
\author[E. Lerman]{Eugene Lerman}

\address{Inst.\ of Mathematics, The Hebrew Univ.\ of Jerusalem, 
Giv'at Ram, Jerusalem 91904, Israel}
\address{Dept.\ of Mathematics, Univ.\ of Illinois, Urbana, IL 61801}

\thanks{dg-ga/9506008}

\thanks{The work of E. Lerman was partially supported by an NSF
postdoctoral fellowship. The work of Y. Karshon was partially
supported by NSF grant DMS-9404404.}

\email{karshon@@math.huji.ac.il}
\email{lerman@@math.uiuc.edu}


\date{\today}

\maketitle

\begin{abstract}
Let $\Phi : M \to \g^*$ be a proper moment map associated to an action
of a compact connected Lie group, $G$, on a connected symplectic manifold,
$(M,\omega)$.  A {\em collective function\/} is a pullback via $\Phi$
of a smooth function on $\g^*$. In this paper we present four new
results about the relationship between the collective functions and
the $G$-invariant functions in the Poisson algebra of smooth functions on
$M$.  More specifically, we show:

\noindent
1. \ 
The centralizer of the invariant functions consists of the algebra 
of smooth functions on $M$ that are constant on the level sets of 
the moment map.
This resolves a conjecture of Guillemin and Sternberg. 

\noindent
2. \ The question of whether this centralizer is equal to the algebra 
of collective functions or is larger is equivalent to a formal
algebraic question on the level of power series.

\noindent
3. \ If the group $G$ is a torus, the centralizer of the
invariant functions consists of the collective functions.  We close a
gap in earlier proofs of this fact.

\noindent
4. \ 
If the group  $G$ is $SU(2)$ and the centralizer of the invariant
functions is larger than the algebra of of collective functions,
the action of $SU(2)$ extends to an action of $U(2)$ with the same orbits, 
and the centralizer of the invariant functions consists of the 
$U(2)$-collective functions.
\end{abstract}

\tableofcontents

\section{Introduction}
Let $\Phi : M \to \g^*$ be a moment map associated to a Hamiltonian
action of a compact connected Lie group $G$ on a compact
\footnote{
Throughout this introduction we assume that the manifold $M$
is compact and the group $G$ is connected.
In the rest of the paper our assumptions are often more general.
}
connected symplectic manifold $(M,\omega)$.
Pullbacks by $\Phi$ of smooth functions on $\g^*$ are called
{\em collective functions}. They form a Poisson subalgebra
of the algebra of smooth functions on $M$.  Its centralizer is the
algebra of invariant functions, i.e., a smooth function $f$ on $M$ is
invariant if and only if $\{ f,h \} =0$ for every collective function
$h$, where $\{ , \}$ denotes the Poisson bracket corresponding to the
symplectic form $\omega$.

Motivated by a study of completely integrable systems in
\cite{g-s:collective}, Guillemin and Sternberg conjectured in
\cite{g-s:mfs} that the centralizer of the algebra of invariant functions
is the algebra of collective functions.  They proved this conjecture for
neighborhoods of generic points in $M$.

A collective function is clearly constant on the level sets of the
moment map. The converse need not be true. For example, the standard
linear action of the group $G=SU(2)$ on $\C^2$ has a moment map
$\Phi(u,v) = ( \ol{u} v , \half |u|^2 - \half |v|^2 )$
when we identify the vector
space $\g^*$ with $\R \times \C$.  The function $f(u,v) = |u|^2 +
|v|^2$ is constant on the level sets of $\Phi$ because it is equal to
$ (|\ol{u} v|^2 + (\half |u|^2 - \half |v|^2)^2 )^\half = \half ||\Phi||$.
It is not collective because the function $||x||$ is not smooth 
on $\R \times \C$.

In section \ref{sec:duality} of this paper we show that the
centralizer of the algebra of invariant functions is the algebra of
functions that are constant on the level sets of the moment map.  In
fact, these two algebras are mutual centralizers in the Poisson
algebra $\Cinf(M)$. See Theorem \ref{duality} and Corollary
\ref{mutual-conn}.  This was already shown in the thesis of the first
author \cite{k:thesis}, but our current proof is shorter.

This result raises the following question: what is the obstruction for a
function that is constant on the level sets of the moment map to
be collective?  In section \ref{sec:BM}, Theorem \ref{pullback},
we express this obstruction as a condition on the Taylor series 
of the function. The proof uses theorems of Bierstone and Milman 
and of Marle, Guillemin, and Sternberg.  Theorem \ref{pullback}
essentially reduces the identification of the centralizer of the 
invariant functions to an algebraic question. Based on this,
F. Knop recently announced a complete description of the
centralizer of the invariant functions on a Hamiltonian space in terms
of the little Weyl group (defined in \cite{knop:1990}) of the space.

If the Lie group $G$ is abelian, every function which is constant on
the level sets of the moment map is a collective function.
So the conjecture of Guillemin and Sternberg is true for torus actions.
The history of the proof is as follows.
Already in their paper \cite{g-s:mfs}, Guillemin and Sternberg
proved that for a linear torus action on a symplectic vector space
the centralizer of the invariants consists of
functions that are constant on the level sets of the moment map. 
They claimed that these functions are collective.  
This claim is not obvious; we prove it in section \ref{sec:toral}
of this paper. In \cite{l:centralizer}, the second author 
showed that this claim implies that the conjecture of Guillemin and 
Sternberg is true for torus actions on compact manifolds, and, more
generally, for actions of compact Lie groups on compact manifolds,
provided that the image of the moment map does not intersect the walls
of the Weyl chambers.  We recall (slightly stronger versions of) these
results in section \ref{sec:toral}.

The first counterexample to the conjecture of Guillemin and Sternberg
was given by the second author in \cite{l:centralizer}.  This is the
standard action of $SU(2)$ on $\C^2$ with moment map $\Phi(u,v) =
( \ol{u}v , \half |u|^2 - \half |v|^2 )$. This action
of $SU(2)$ extends to the standard action of $U(2)$ on $\C^2$ with
the same orbits (spheres) and with a moment map $\tPhi(u,v) = 
( \ol{u} v , |u|^2 , |v|^2 )$.  The centralizer of the invariants 
(for either $SU(2)$ or $U(2)$) consists of the $U(2)$-collective functions.
For instance, the function $f(u,v) = |u|^2 + |v|^2$ on $\C^2$
is in the centralizer of the invariants, is not $SU(2)$-collective, 
but is $U(2)$-collective.  

A similar phenomenon happens more generally:
for a Hamiltonian action of $SU(2)$ with a proper
moment map, either the centralizer of the invariant functions
consists of the collective functions, or the action of $SU(2)$
extends to an action of $U(2)$ with the same orbits and for which the
centralizer of the invariant functions consists of the $U(2)$-collective
functions. This we show in section \ref{sec:SU2}.

In the rest of this introduction we describe the context in which
Guillemin and Sternberg posed their conjecture.
The notion of mutually centralizing subgroups in the symplectic group
originated in physics. It was studied by Sternberg and Wolf, by Howe,
by Kashiwara and Vergne, and by Jakobsen and Vergne, in \cite{s-w},
\cite{howe}, \cite{kash-vergne}, and \cite{jak-vergne} respectively.
In the classical analogue of this notion one considers two connected
Lie groups, $G$ and $H$, that act on a symplectic manifold $M$ with
moment maps
\begin{equation} \labell{eq2}
\begin{array}{c}
 M \\
 \stackrel{F}{\swarrow} \phantom{M} \stackrel{\Phi}{\searrow} \\
 \h^* \phantom{MMM} \g^*
\end{array}
\end{equation}
such that the respective algebras of collective functions, 
$F^* C^\infty(\h^*)$ and $\Phi^* C^\infty(\g^*)$, 
are mutual centralizers in the Poisson algebra $C^\infty (M)$.  
The $G$-moment map, $\Phi$, then becomes an orbit map
for the action of $H$. This means that the $H$-invariants are exactly
the pull-backs via $\Phi$ of smooth functions on $\g^*$.  Similarly,
the $G$-invariants are the $H$-collective functions.  
This has several consequences. First, 
generically, the $G$-reduced spaces are coverings of coadjoint orbits
of $H$ and vice versa.  Moreover, we get a correspondence between the
coadjoint orbits of $G$ and those of $H$ which occur in the images of
the respective moment maps.
This phenomenon, which is the classical analogue of Howe's dual pairs
in representation theory \cite{howe}, was observed and explained by
Kazhdan, Kostant and Sternberg in \cite{k-k-s} for the case that the
$G$-orbits form a foliation.  

More generally, A. Weinstein \cite{w:poisson} defined a {\em dual
pair\/} to be a pair of Poisson maps $f: M \to A$ and $g: M \to B$
from a symplectic manifold $M$ to Poisson manifolds $A$ and $B$
such that the algebras $f^* \Cinf(A)$ and $g^* \Cinf(B)$ are mutual
centralizers in the Poisson algebra $\Cinf(M)$.  Dual pairs and their
infinite dimensional analogs occur in the study of tops, compressible
fluids, elasticity, Maxwell - Vlasov equations \cite{MRW}, etc., 
and have lead to the notion of Morita equivalence of Poisson manifolds
(cf.\ \cite{g-l}).

The conjecture of Guillemin and Sternberg is equivalent to the maps
$$ \begin{array}{c} 
   M \\ 
   \stackrel{\pi}{\swarrow} \phantom{M} \stackrel{\Phi}{\searrow} \\
   M/G \phantom{MMM} \g^*
   \end{array} 
$$
forming a ``dual pair'', where $\Phi$ is the moment map and $\pi$ is the
quotient map, and where we interpret the Poisson algebra $\Cinf(M/G)$
as the algebra of functions on $M/G$ whose pullback to $M$ is
smooth. Note that the quotient $M/G$ need {\em not} be a manifold
hence the quotes around the expression ``dual pair.''
 
Given a Hamiltonian action of a Lie group $G$, one may wonder whether
there exists a Hamiltonian action of another Lie group, $H$, such that
the corresponding moment maps form a dual pair \eqref{eq2}. 
A necessary condition for
the existence of this other action is that the centralizer in $\Cinf(M)$
of the algebra of $G$-invariant functions be equal to the algebra of 
$G$-collective functions.

An example to keep in mind is
the standard action of the group $G=U(2)$ on the symplectic vector
space $M=\C^2$. The $G$-orbit map is $F(z,w) = |z|^2 + |w|^2$, which
generates the diagonal action of $H=S^1$.
Another interesting example is the natural action of the orthogonal
group $O(k)$ on $(T^* {\Bbb R}^k)^n$, the $n$-fold product of the
cotangent bundle of ${\Bbb R}^k$.  The space $(T^* {\Bbb R}^k)^n$ is
the phase space of the $n$-body problem, and $O(k)$ is its natural
symmetry group.  The group $H$ in this case is the symplectic group
$Sp(\R^{2n})$ (cf.\ \cite{LMS}).

\section{The centralizer of invariant functions}
\labell{sec:duality}
Let $\Phi : M \to \g^*$ be a moment map associated to a Hamiltonian
action of a compact Lie group $G$ on a symplectic manifold $(M, \omega)$.
Recall, this means that for any element $\xi$ of the Lie algebra $\g$ 
of $G$ we have $d \Phi^\xi = - \iota (\xi_M) \omega$,
where $\Phi^\xi = \langle \Phi, \xi \rangle$
is the $\xi$-component of the moment map 
and $\xi_M$ is the vector field on $M$ that generates the action
of the one parameter subgroup $\{\exp(t \xi) \, , \, t \in \R\}$ of $G$.
We also require that $\Phi$ be equivariant with respect to 
the given action of $G$ on $M$ and the coadjoint action on $\g^*$.
The main result of this section reads:

\begin{Theorem} \labell{duality}
The centralizer of the algebra of $G$-invariant functions 
in the Poisson algebra of smooth functions on $M$
is the set of smooth functions that are locally constant 
on the level sets of the moment map.
\end{Theorem}

\begin{pf}
Since the Hamiltonian flow of an invariant function preserves the level
sets of the moment map, the Poisson bracket of an invariant function and
a function that is locally constant on the level sets of the moment map
is zero.  This shows that the centralizer of the invariant functions
contains the functions that are locally constant on the level sets of
the moment map.  We would like to show that there is nothing else in
the centralizer.

Let $h$ be a function in the centralizer of the invariant functions.
Let $\gamma(t)$ be any smooth curve contained in a level set of
the moment map $\Phi$.  Since any two points in a connected component
of a level set of $\Phi$ can be connected by a piece-wise smooth curve
(see Lemma~\ref{smooth-curve}), we are done if we can prove that the
derivative of $h(\gamma (t))$ is zero for all $t$.  This derivative
is equal to $\omega(\dot{\gamma},X_h)$ where $X_h$ is the Hamiltonian
vector field of $h$.

For any vector $\xi$ in the Lie algebra $\g$ we have
$0 = \langle \dot{\gamma},d\Phi^\xi \rangle = \omega(\dot{\gamma},\xi_M)$.
Hence if $\gamma (t)$ is a smooth curve contained in a level set
of the moment map, the tangent vectors $\dot{\gamma}$ lie
in the symplectic perpendiculars to the $G$-orbits.

To finish the argument it suffices to show that the Hamiltonian vector
field, $X_h$, of a function $h$ in the centralizer of the invariant
functions is tangent to the $G$-orbits.  Let $\sigma(t)$ be an integral
curve of the vector field $X_h$.  Then for any $G$-invariant function,
$f$, we have $\frac{d}{dt} (f(\sigma (t)) = (X_h f) (\sigma (t)) = 0$,
i.e., $f$ is constant along $\sigma(t)$.  Since, $G$ being compact, the
$G$-invariant functions separate orbits, the integral curve $\sigma (t)$
is contained in a single $G$-orbit. Hence the vector field $X_h$ is
tangent to $G$-orbits.  This proves Theorem~\ref{duality}.  \end{pf}

The rest of this section contains corollaries of Theorem \ref{duality}.
We assume throughout that $\Phi$ is the moment map for a Hamiltonian
action of a compact Lie group $G$ on a symplectic manifold $(M,\omega)$.
We do not assume that $G$ is connected or that $\Phi$ is proper
unless we explicitly say so.

\begin{Corollary}
The set of functions that are locally constant on the
level sets of $\Phi$ is a Poisson algebra.
\end{Corollary}

\begin{pf}
The Jacobi identity implies that for any subset, $R$, of $\Cinf(M)$
the centralizer, $R^c := \{ f \in\Cinf(M) \ | \ \{ f,h \} =0
\mbox{ for all } h \in R \}$, is a Poisson algebra.  
Apply this when $R$ is the set of invariant functions.
\end{pf}

The conjecture of Guillemin and Sternberg in \cite{g-s:mfs} is stated
differently than the way we quoted it in the introduction.  Namely, it
is stated in terms of the the double centralizer (``double commutator")
of the set of functions $\{ \Phi_1, \ldots, \Phi_{\dim G} \}$, where
$\Phi_i$ are the coordinates of the moment map with respect to some
basis of the vector space $\g^*$.  However,

\begin{Lemma} \labell{magavot}
The double centralizer of the set of coordinate functions of the moment
map is equal to the centralizer of the set of invariant functions.
\end{Lemma}

The proof uses two lemmas:

\begin{Lemma} \labell{easy}
A function Poisson commutes with all the coordinates of the moment map if
and only if the function is invariant under the identity component of the
group.  
Consequently, the centralizer of the collective functions is equal
to the functions invariant under the identity component of the group.
\end{Lemma}

\begin{pf}
This follows easily from the definition of the moment map.
\end{pf}

\begin{Lemma} \labell{G0}
Every smooth function on $M$ which is invariant under the identity
component of the compact group $G$ is locally equal to a $G$-invariant
function.
\end{Lemma}

\begin{pf}
This is an easy consequence of the slice theorem for actions of compact 
groups.
\end{pf}

\begin{pf*}{Proof of Lemma \ref{magavot}} 
If $G$ is connected, the lemma follows immediately from Lemma \ref{easy}.
If $G$ is disconnected, the centralizer of the $G$-invariant functions 
is the same as the centralizer of the $G_0$-invariant functions
where $G_0$ is the identity component; this follows from Lemma \ref{G0}.
\end{pf*} 

\begin{Notation} \labell{notation} 
Let $\Cinf(M)^G$ denote the $G$-invariant smooth functions on $M$.
Let $G_0$ be the identity component of $G$, and
let $\Cinf(M)^{G_0}$ denote the $G_0$-invariant smooth functions;
equivalently, these are the smooth functions on $M$ that are locally constant
on the $G$-orbits.

Let $\Cinf(M)^\Phi$ denote the smooth functions on $M$ that are
constant on the level sets of the moment map $\Phi$, and
$\Cinf(M)^\Phi_\loc$ the functions that are {\em locally\/} constant
on these level sets.  Let $\Phi^* \Cinf(\g^*)$ denote the collective
functions, i.e., pullbacks by $\Phi$ of smooth functions on $\g^*$.
\end{Notation}

In this notation, Theorem~\ref{duality} says that
the centralizer of $\Cinf(M)^G$ is $\Cinf(M)^\Phi_\loc$.

\begin{Corollary} \labell{mutual}
The algebra $\Cinf(M)^\Phi_\loc$ and the algebra $\Cinf(M)^{G_0}$
are mutual centralizers in the Poisson algebra $\Cinf(M)$.

Moreover, the centralizers of the algebras
$\Phi^* \Cinf(M) \subseteq \Cinf(M)^\Phi \subseteq \Cinf(M)^\Phi_\loc$
are all equal to $\Cinf(M)^{G_0}$,
and the centralizers of the algebras $\Cinf(M)^G \subseteq \Cinf(M)^{G_0}$
are both equal to $\Cinf(M)^\Phi_\loc$.
\end{Corollary}

\begin{pf}
Theorem \ref{duality} applied to the group $G_0$ implies that the 
centralizer of $\Cinf(M)^{G_0}$ is equal to $\Cinf(M) ^\Phi_\loc$.
Conversely, since the algebra $\Cinf(M)^\Phi_\loc$ contains the collective
functions, its centralizer is contained in the centralizer of the
collective functions, which is equal to the $G_0$-invariant functions
by Lemma~\ref{easy}.
\end{pf}

\begin{Corollary} \labell{dual}
Let $\Phi$ be a moment map associated to a Hamiltonian action of a
compact Lie group $G$ on a symplectic manifold $(M,\omega)$.
If $G$ is connected, the following properties of $\Phi$ are equivalent
to each other:
\begin{enumerate}
\item
The algebra of collective functions and the algebra of invariant functions 
are mutual centralizers in the Poisson algebra $\Cinf(M)$
\item
The double centralizer of the set of coordinate functions
of the moment map is the algebra of collective functions
\item
Every smooth function on $M$ that is locally constant on the level sets 
of the moment map is collective.
\end{enumerate}
\end{Corollary}

\begin{pf}
Since the centralizer of the collective functions is the invariant 
functions (Lemma \ref{easy}), condition 1 is equivalent to
\begin{enumerate}
\item[1'.]
{\em
The centralizer of the algebra of invariant functions
is the algebra of collective functions.
}
\end{enumerate}
The equivalence of conditions 1' and 2 is immediate from Lemma
\ref{magavot}; the equivalence of conditions 1' and 3 is immediate from
Theorem \ref{duality}.
\end{pf}

\begin{Definition} \labell{property} 
A smooth map $\Phi: M \to N$ between two smooth manifolds
has the {\em division property\/} if any smooth function on $M$ that
is locally constant on the level sets of $\Phi$ is the pullback via
$\Phi$ of a smooth function on $N$.
\end{Definition}

The hard part of Corollary \ref{dual} can be rephrased as follows:

\begin{Corollary} \labell{dual2}
Let $\Phi : M \to \g^*$ be a moment map associated to an action of a
compact connected Lie group $G$ on a symplectic manifold
$(M,\omega)$. The algebra of collective functions and the algebra of
invariant functions are mutual centralizers in the Poisson algebra
$\Cinf(M)$ {\em if and only if\/} the moment map $\Phi$ has the
division property.
\end{Corollary}

One obstruction for a map to have the division property is topological
--- the connectedness of the level sets of the map.  Another obstruction
has to do with analytic properties of the map.  The following examples
illustrate these ideas: (1) The map $x \mapsto x^3$ on $\R$ does not
have the division property because of the singularity at $x=0$; the
function $f(x)=x$ is not a pullback of a smooth function.  (2) The map
$e^{i\theta} \mapsto e^{2i\theta}$ on the unit circle does not have the
division property because the level sets are not connected; the function
$f(e^{i\theta})=\cos(\theta)$ is not a pullback.  (3) The map $(x,y)
\mapsto x^2 + y^2$ from $\R^2$ to $\R$ has the division property; every
rotationally invariant smooth function is a smooth function of $x^2 +
y^2$ . This is a special case of a theorem of G. Schwarz \cite{sch:sm}.

Two global properties of a moment map will be relevant to us: 
properness of the map and connectedness of its level sets.  
The first of these often implies the second:

\begin{Proposition} \labell{connected0}
Let $\Phi : M \to \g^*$ be a moment map associated to an action of a
compact Lie group, $G$, on a connected symplectic manifold, $(M,\omega)$.
If this map is proper, its level sets are connected.
\end{Proposition}

\begin{pf}
For torus actions this proposition was proved by Atiyah  \cite{at}.
The general case was proved by Kirwan \cite{Kirwan,Kirwan:book}.
\end{pf}

This result has recently been generalized:
\begin{Proposition} \labell{connected}
Let $\Phi : M \to \g^*$ be a moment map associated to an action of a
compact Lie group $G$ on a connected symplectic orbifold $(M,\omega)$.
Let $\t^*_+$ be a Weyl chamber, identified with a subset of $\g^*$ via a
choice of an ${\text Ad}$-invariant metric.  Suppose that there exists
a $G$-invariant open subset $N$ of $\g^*$ containing the image $\Phi
(M)$, such that the intersection $N \cap \t^*$ is convex and the map
$\Phi$ is proper as a map from $M$ to $N$. Then the level sets of $\Phi$
are connected.
\end{Proposition}

\begin{pf}
See \cite{LMTW}.
\end{pf}

We proceed with consequences of Theorem \ref{duality}.

\begin{Corollary} \labell{mutual-conn}
Let $\Phi : M \to \g^*$ be a moment map associated to a Hamiltonian
action of a compact Lie group, $G$, on a connected symplectic manifold,
$(M,\omega)$.  If the group $G$ is connected and the moment map, $\Phi$,
is proper, the algebras $\Cinf(M)^\Phi$ and $\Cinf(M)^G$ are mutual
centralizers in the Poisson algebra $\Cinf(M)$.
(See Notation \ref{notation}.)
\end{Corollary}

\begin{pf}
If $G$ is connected, $\Cinf(M)^{G_0} = \Cinf(M)^G$.  If $\Phi$ is
proper, its level sets are connected, so $\Cinf(M)^\Phi_\loc =
\Cinf(M)^\Phi$.  The rest is immediate from Corollary~\ref{mutual}.
\end{pf}

\begin{Corollary} \labell{cor1}
If the moment map, $\Phi : M \to \g^*$, is proper, the centralizer of
the invariant functions consists of those smooth functions on $M$ that
are pullbacks by $\Phi$ of {\em continuous} functions on $\g^*$:
\begin{equation} \labell{c0}
	 (\Cinf(M)^G)^c = \Cinf(M) \cap \Phi ^* C^0 (\g^*).
\end{equation}

Let $N \subset \g^*$ be an open set containing the moment image,
$\Phi(M)$.  If the moment map, $\Phi$, is proper as a map from $M$ to
$N$ and has connected level sets, the centralizer of the invariant
functions consists of those smooth functions on $M$ that are pullbacks
by $\Phi$ of continuous functions {\em on $N$}:
$$ 	
(\Cinf(M)^G)^c = \Cinf(M) \cap \Phi ^* C^0 (N).
$$
\end{Corollary}

\begin{pf*}{Proof of Corollary \ref{cor1}}
The first part follows from the second, because the level sets of a proper
moment map are connected (Proposition \ref{connected0}).  A function in
the centralizer of the invariants is the pullback of a function on $N$
(Theorem~\ref{duality}).  The function on $N$ is continuous because the
moment map is proper.
\end{pf*}

\begin{Remark}
If the moment map is not proper, \eqref{c0} may fail to hold. 
See Examples \ref{counter1} and \ref{counter2}.  
\end{Remark}

\section{Division property can be detected formally}
\labell{sec:BM}

Let $\Phi : M \to \g^*$ be a proper moment map associated
to an action of a compact Lie group $G$ on a symplectic 
manifold $(M,\omega)$. In this section we show that the difference
between the algebra of functions that are constant on the level sets
of the moment map and the algebra of collective functions
can already be detected on the level of power series of these functions.
This result is stated in Theorem~\ref{pullback}. The main idea of the
proof is to apply a theorem of Bierstone and Milman to the 
Marle-Guillemin-Sternberg local normal form of the moment map.

We first recall a definition of Bierstone and Milman.  Recall that a
function is {\em flat\/} at a point if its Taylor series at that point
vanishes.

\begin{Definition}  \labell{def:formal-pullback}
Let $\Phi: M \to N$ be a smooth map between two smooth manifolds.
A smooth function $f$ on $M$ is {\em a formal pullback with respect
to $\Phi$\/} 
if for every point $y$ in the image $\Phi(M)$ there exists a
function, $\varphi$, on $N$ such that 
$f - \Phi^* \varphi$ is flat at all the points of $\Phi\inv(y)$.
\end{Definition}

\begin{Remark} \labell{constant}
Every formal pullback with respect to $\Phi$ is constant
on the level sets of $\Phi$; if $f - \Phi^* \varphi$ is 
flat, $f(x) = \varphi(y)$ for all $x \in \Phi\inv(y)$.
\end{Remark}

\begin{Remark}
In the notation of Bierstone and Milman, the set of formal pullbacks
with respect to $\Phi$ is $(\Phi^*\Cinf(N))\hat{}$; they have no term
to describe the elements of this algebra.  In a more recent paper
\cite{BMP}, with W. Pawlucki, they use the term ``formal composite
with $\Phi$" to describe a formal pullback.
\end{Remark}

\begin{Remark}
Recall, the Taylor series at a point $p$ of a smooth function 
$f : M \to \R$ is an element of the algebra $\prod_i S^i(T^*_pM)$ of
(formal) power series at $p$. Conversely, by Borel's theorem, every
power series at $p$ is a Taylor series of some smooth function on $M$.
Therefore, one can think of a power series at a point $p$ as an
equivalence class of functions:  two functions are equivalent if and
only if their difference is flat at $p$.  Since the pullback of functions 
induces a well defined pullback of Taylor series, being a formal pullback
is a condition on Taylor series:
a smooth function $f : M \to \R$ is a formal pullback with respect
to a smooth map
$\Phi : M \to N$ if and only if for every $y \in N$ there exists a power 
series $\varphi$ on $N$, centered at $y$, such that for all $x$ in the level
set $\Phi\inv(y)$, the power series of $f$ at $x$ is the pullback
of the power series $\varphi$.
\end{Remark}

We will now state the main result of this section.  Recall, a continuous
map $\psi :A \to B$ is {\em semi-proper} if for every compact set
$L\subset B$ there is a compact set $K \subset A$ such that $\psi (K)
= L\cap \psi (A)$.

\begin{Theorem}  \labell{pullback}
Let $\Phi : M \to \g^*$ be a moment map associated to a Hamiltonian
action of a compact Lie group $G$ on a connected symplectic manifold
$(M,\omega)$.  If this map, $\Phi$, is proper, every formal pullback
with respect to $\Phi$ 
is a collective function, i.e., is in $\Phi^* \Cinf(\g^*)$.

Let $N$ be an open subset of $\g^*$ containing the moment image,
$\Phi(M)$. If the moment map, $\Phi$, is semi-proper as a map from $M$
to $N$ and has connected level sets, every formal pullback with 
respect to $\Phi$
is a pullback of a smooth function on $N$, i.e., is in $\Phi^* \Cinf(N)$.
\end{Theorem}

\begin{Corollary} \labell{criterion}
A proper moment map $\Phi$ has the division property
(see Definition \ref{property})
if and only if every smooth function on $M$ that is locally 
constant on the level sets of $\Phi$ is a formal pullback with 
respect to $\Phi$.
\end{Corollary}

The proof of Theorem~\ref{pullback} relies upon the following
theorem of Bierstone and Milman.  Recall that a {\em semi-analytic}
subset of an analytic manifold is a subset that is locally defined by
inequalities involving analytic functions.

\begin{Theorem}[Bierstone-Milman {\cite[Theorem 0.1]{BM2}}] 
	\labell{BM} 
Let $M$ and $N$ be real analytic manifolds. Let $\Phi: M \to N$ be a
real analytic mapping that is semi-proper and whose image, $\Phi(M)$,
is semi-analytic.  Then a function $f$ is a formal pullback with 
respect to $\Phi$
if and only if it is the pullback by $\Phi$ of a smooth function on $N$.
\end{Theorem}

\begin{Remark} 
Theorem~0.1 in \cite{BM2} requires the image $\Phi(M)$ to be 
Nash subanalytic.  Bierstone and Milman point out in \cite{BM1}
that every semi-analytic set is Nash subanalytic.
\end{Remark}

A priori, our manifolds and maps are only smooth and not real analytic,
so we cannot apply Theorem~\ref{BM} directly. We will use the following
variant:

\begin{Proposition} \labell{loc-BM}
Let $M$ and $N$ be smooth manifolds, and let $\Phi:M \to N$ be a smooth 
map that satisfies the following conditions.
\begin{enumerate}
\item	The image $\Phi(M)$ is closed.
\item
	For every point $x$ in $M$ there exist neighborhoods $U_x$
	of $x$ in $M$ and $W_x$ of $\Phi(x)$ in $N$ such that
\begin{enumerate}
	\item	$\Phi(U_x) = \Phi(M) \cap W_x$;
	\item	the restriction $\Phi |_{U_x} : U_x \to W_x$
		is semi-proper; 
	\item	there exist real analytic structures on $U_x$ and on $W_x$
		compatible with their smooth structures such that the 
		restriction $\Phi|_{U_x} : U_x \to W_x$ is a real analytic
		map whose image is a semi-analytic subset of $W_x$.
\end{enumerate}
%\item The level sets of the map $\Phi$ are connected.
\end{enumerate}
Then the set of pullbacks by the map $\Phi$ 
coincides with the set of formal pullbacks with respect to 
$\Phi$. 
\end{Proposition}

\begin{pf}
Clearly, every pullback is a formal pullback. Conversely, let $f \in
\Cinf (M)$ be a formal pullback with respect to $\Phi$.  Let $x$ be a
point in $M$, and let $U_x$ and $W_x$ be as in Condition 2 above.
 Since $f$ is a formal pullback with respect to
$\Phi$, its restriction $f |_ {U_x}$ is a formal pullback with respect
to the map $\Phi|_{U_x} : U_x \to W_x$.  Theorem \ref{BM} applies to
this map because of Conditions 2(b) and 2(c). Hence there exists a
smooth function $\varphi _x$ on $W_x$ such that $f = \varphi_x \circ
\Phi$ on $U_x$.  This equality holds on all of $\Phi\inv(\Phi(U_x))$
because $f$, being a formal pullback with respect to $\Phi$, is
constant on the level sets of $\Phi$ (see Remark \ref{constant}).
Condition 2(a) implies that $\Phi\inv(\Phi(U_x)) = \Phi\inv(W_x)$, so
$f = \varphi_x \circ \Phi$ on all of $\Phi^{-1}(W_x)$.  The open sets
$W_x$ together with the complement of the image $\Phi(M)$ form an open
cover of the target manifold, $N$.  Using a partition of unity
subordinate to this cover we piece together the functions $\varphi_x$
to form a function $\varphi$ on $N$ such that $f = \varphi \circ
\Phi$.
\end{pf}

Thus to prove Theorem \ref{pullback} it is enough to verify that 
every proper moment map on a connected symplectic manifold
satisfies the hypotheses of Proposition \ref{loc-BM}.

\begin{Proposition} \labell{properties}
Let $\Phi : M \to \g^*$ be a moment map associated to an action of a
compact connected Lie group $G$ on a connected
symplectic manifold $(M,\omega)$.
Assume that the map $\Phi$ is semi-proper as a map into an open subset 
$N$ of $\g^*$ and that its level sets are connected.
Let $x$ be a point in $M$, and let $G_\alpha$ be the stabilizer of its
image, $\alpha = \Phi(x)$, under the coadjoint action.  Then there exist a
neighborhood $U_x$ of the orbit $G_\alpha \cdot x$ in $M$ and a
neighborhood $W_x$ of the point $\Phi(x)$ in $\g^*$ with the following
properties.
\begin{enumerate}
\item 	$\Phi(U_x) = \Phi(M) \cap W_x$.
\item	The restriction $\Phi|_{U_x} : U_x \to W_x$
	is semi-proper.
\item	There exist real analytic structures on $U_x$ and on $W_x$,
	compatible with their smooth structures, such that
	\begin{enumerate}
	\item	the restriction $\Phi|_{U_x} : U_x \to W_x$
		is a real analytic map;
	\item	the image $\Phi(U_x)$ is a semi-analytic subset of $W_x$.
	\end{enumerate}
\end{enumerate}
Moreover, the neighborhoods $U_x$ and $W_x$ can be chosen to be
arbitrarily small, i.e., can be chosen to be contained in any given
neighborhoods $U'$ of $G_\alpha \cdot x$ and $W'$ of $\Phi(x)$.
\end{Proposition}

Note that only property 1 is global.  It is only to prove this property 
that we assume that the moment map is semi-proper and that its level sets 
are connected.

\begin{pf*}{Proof of Proposition \ref{properties}}
Let us first prove properties 1--3 when the orbit $G \cdot x$ is
isotropic, equivalently, when $\alpha$ is fixed under the coadjoint
action of $G$.  Since $\alpha$ is fixed, the translation $\Phi
-\alpha$ of the moment map by $-\alpha$ is still a moment map.
So, without loss of generality, we can assume that $\alpha = 0$.

In the appendix (Theorem~\ref{locsymp}) we describe a local model
for a neighborhood of an isotropic orbit $G\cdot x$, 
$$ 
	Y = G \times_{G_x} (\g_x^0 \oplus V),
$$
where $G_x$ is the stabilizer of $x$, $\g_x$ is its Lie algebra,
$\g_x^0$ is the annihilator of $\g_x$ in $\g^*$, and the vector space
$V$ is the symplectic slice at $x$.  The action of $G$ on $Y$ is
Hamiltonian with a moment map $\Phi_Y : Y \to \g^*$ given by the
formula
$$ \Phi_Y([g,\eta,v]) = 
   {\text Ad}^\dagger(g) (\eta + i(\Phi_V(v))),$$
where $\Phi_V$ is a quadratic map from $V$ to $\g_x^*$
and $i$ is a $G_x$-equivariant embedding of $\g_x^*$ in $\g^*$.
Moreover, by Theorem~\ref{locsymp} there exists a neighborhood $U_x$ 
of $G \cdot x$ in $M$ and an equivariant embedding, 
$\iota : U_x \to Y$, of $U_x$ onto a
neighborhood of the zero section in the  model $Y$, such that
$\Phi = \Phi_Y \circ \iota$.

By Lemma \ref{local-properties}, the image under the moment map of a
small neighborhood of an orbit $G \cdot x$ does not change as $x$
varies along a connected component of the level set
$\Phi\inv(\Phi(x))$.  By Remark~\ref{cone}, this image is the
intersection of the cone $\Phi_Y(Y)$ with a neighborhood of the origin
in $\g^*$.  Lemma \ref{local-properties} together with the facts that
$\Phi$ is semi-proper and that its level sets are connected implies
that we can choose a neighborhood $W_x$ of the origin in $\g^*$ and
shrink the neighborhood $U_x$ of $G\cdot x$ so that $\Phi(U_x) =
\Phi(M) \cap W_x = \Phi_Y(Y) \cap W_x$, i.e., so that property 1
holds.

The map $\Phi_Y$ is analytic with respect to the natural
real analytic structures of the model $Y$ and of the vector space $\g^*$.
If we endow $U_x$ with the real analytic structure induced by its 
embedding, $\iota$, into $Y$, property 3(a) holds.

Consider the action of $\R_+$ on $Y$ given by
$$	
\lambda \cdot [g,\eta,v] = [g, \lambda \eta, \sqrt{\lambda} v].
$$
The map $\Phi_Y : Y \to \g^*$ is homogeneous of degree one with
respect to this action of $\R_+$.  After possibly shrinking $U_x$ and
$W_x$ further, we can assume that the open set $\iota(U_x) \subseteq
Y$ is preserved under multiplication by any $\lambda < 1$; for
such $\lambda$ we define $\lambda : U_x \to U_x$ by $\iota(\lambda
\cdot m) = \lambda \cdot \iota(m)$.  Let $K$ be a compact subset of
the open set $W_x$.  Then there exists a positive number $\lambda < 1$
such that $K$ is contained in $\lambda W_x$.  By homogeneity, $K \cap
\Phi(U_x)$ is contained in $\Phi(\lambda \cdot U_x).$ Then $L :=
\text{closure}(\lambda \cdot U_x) \cap \Phi\inv(K)$ is a compact
subset of $U_x$ whose image is $K \cap \Phi(U_x)$.  This proves
property 2.

Since the map $\Phi_V$ is algebraic, its image, $\Phi_V(V)$, is a
semi-algebraic subset of $\g_x^*$, by the Tarski-Seidenberg theorem
(see, for example, \cite[Theorem~2.3.4]{BR}).  Furthermore, since
$\text{Ad}^\dagger(G) \subseteq \text{GL}(\g^*)$ is algebraic, the set
$\Phi_Y(Y) = \text{ Ad}^\dagger(G)(\g_x^0 \times \Phi_V(V))$ is a
semi-algebraic subset of $\g^*$.  Restricting to the open subset
$W_x$, we see that $\Phi(U_x) = \Phi_Y(Y) \cap W_x$ is a semi-analytic
subset of $W_x$. This proves property 3(b).

We now remove the assumption that the orbit is isotropic.
Let $\alpha = \Phi (x)$, and
let $G_\alpha$ denote its stabilizer under the coadjoint action. Let $S$
be a slice at $\alpha$ for the action of $G$ and $R= \Phi \inv (S)$ the
corresponding symplectic cross-section (cf.\ Theorem~\ref{cross-section}).
By Corollary~\ref{use_cross-section}, up to a composition with 
diffeomorphisms, the moment map $\Phi$ is a fiber bundle map
$$
  G\times _{G_\alpha} R \to G\times _{G_\alpha} S, \quad\quad [g, r] 
  \mapsto [g, \Phi_R (r)], 
$$ 
and $\Phi_R = \Phi |_R$ is the $G_\alpha$ moment map.

The subgroup $G_\alpha$ acts on $G$ by right multiplication.  Since $G
\to G/G_\alpha$ is a locally trivial fibration, there exists a section
on a neighborhood $\calV$ of the identity coset in $G / G_\alpha$.
This section simultaneously trivializes the bundles $G\times
_{G_\alpha} R \to G/G_\alpha$ and $G\times _{G_\alpha} S \to
G/G_\alpha$ over the set $\calV$.
With respect to these trivializations, the moment map 
$\Phi$ is the map
$$      
\text{id} \times \Phi_R: \calV \times R \to \calV \times S\quad\quad 
	(n,r) \mapsto (n,\Phi_R(r)).
$$
Since $G_\alpha \cdot x$ is isotropic in $R$ (Remark~\ref{iso_remark})
and $\Phi_R$ is a $G_\alpha$ moment map, properties 1--3 are satisfied
by the map $\Phi_R : U_x \to W_x$ where $U_x$ is a neighborhood of 
$G_\alpha \cdot x$ in $R$ and $W_x$ is a neighborhood of $\alpha$ in $S$.  
Properties 1--3 for $\Phi_R$ immediately imply properties 1--3
for $\text{id} \times \Phi_R \ : \ \calV
\times U_x \to \calV \times W_x$. Hence the moment
map $\Phi$ has properties 1--3.
\end{pf*}

To prove our main result, Theorem~\ref{pullback}, it suffices to verify
that the hypotheses of Proposition~\ref{loc-BM} are satisfied by the
moment map:

\begin{pf*} {Proof of Theorem~\ref{pullback}}
The first part of the theorem follows from the second by setting $N = \g^*$.

Let $N$ be an open subset of $\g^*$ containing the moment image,
$\Phi(M)$, with the property the moment map $\Phi : M \to N$ is
semi-proper. 
Condition 1 in Proposition \ref{loc-BM}, that the image $\Phi(M)$
is a closed subset of $N$, is satisfied because the image of any
semi-proper map is closed.  Conditions 2(a) -- 2(c) hold by
Proposition~\ref{properties}.  Theorem~\ref{pullback} then follows from
Proposition~\ref{loc-BM}.
\end{pf*}

We conclude this section with two examples which show that 
the properness condition in Theorem~\ref{pullback} is necessary.

\begin{Example} \labell{counter1}
Let $M$ be the cotangent bundle of the two dimensional torus, $T^2$,
minus the zero section. This manifold is the product $M= T^2 \times
(\R^2 \ssminus \{0\})$.  The moment map $\Phi$ is the projection onto
the second factor; it is not proper.  The function $f(x,y) = y /
\sqrt{x^2 + y^2}$ does not extend to a smooth function on $\R^2$ but
it does pull back to a smooth function on $M$ which is a formal
pullback.
\end{Example}

\begin{Example} \labell{counter2}
In this example, a formal pullback $f$ does not even descend
to a continuous function on the image of the moment map.

We construct a Hamiltonian $T^2$-space by gluing two spaces.
The first space, $M_1$, is the product $T^2 \times U$ where $U$
is the subset
of $\R^2$ obtained by removing the origin and the positive $x$-axis.
We can view $M_1$ as an open subset of the cotangent bundle of $T^2$
and take the induced symplectic form.
The moment map is the obvious projection onto $U$.  Its image
is $\R^2$ minus the origin and the  positive $x$-axis.

The second space, $M_2$, consists of the points of $\C^2$ whose first
coordinate is nonzero. This space inherits a symplectic form and an
action of $T^2$ from $\C^2$.  The moment map sends $(z,w) \mapsto
(|z|^2, |w|^2)$.  The image of the moment map is the set
$\{(x, y)\in \R^2: x>0,\,y\geq 0\}$.

We glue the two spaces along the pre-images of the open positive
quadrant by sending $(z,w)$ to
$(\frac{z}{|z|},\frac{w}{|w|}, |z|^2, |w|^2)$.
This gluing map is an equivariant symplectomorphism.

We obtain a space $M$ with a symplectic form, a $T^2$-action and a
moment map. The image of the moment map is $\R^2$ minus the
origin. The branch of $\operatorname{arctan}(y/x)$ which is
discontinuous along the positive $x$-axis pulls back to a smooth
function on $M$.
\end{Example}

\section{Division property of a toral moment map}
\labell{sec:toral}
Consider a Hamiltonian action of a torus, $T$, on a symplectic manifold,
$(M,\omega)$.
Recall that the conjecture of Guillemin and Sternberg, 
asserting that the algebra of invariant
functions and the algebra of collective functions are mutual
centralizers in the Poisson algebra $\Cinf(M)$, holds if and only if
the moment map, $\Phi : M \to \t^*$, has the division property
(Definition \ref{property} and Corollary \ref{dual2}).

Theorem~\ref{pullback}, proved in section~\ref{sec:BM}, provides a
criterion for determining that a moment map has the division property;
see Corollary \ref{criterion}.  In this section we use this criterion
to prove that moment maps arising from torus actions have the division
property.

The results presented here were believed to be known for some time.  
By (the easy part of) Corollary \ref{dual}, our Proposition \ref{torus} 
is equivalent to Proposition~4.1 in \cite{g-s:mfs}, which asserts that
for a symplectic linear action of a torus on a symplectic vector space,
the double centralizer of the set of coordinate functions of the moment 
map is the algebra of collective functions.  Unfortunately, the arguments
presented in \cite{g-s:mfs} only show that this double centralizer
consists of functions that are constant on the level sets of the
moment map.  The missing arguments are not trivial (we don't see a way
of avoiding the theorem of Bierstone and Milman); we provide them in
our proof of Proposition~\ref{torus}.

In \cite{l:centralizer} Lerman used Proposition~4.1 of \cite{g-s:mfs} and
the Marle-Guillemin-Sternberg normal form to deduce  a set of sufficient
conditions for Guillemin-Sternberg conjecture to hold; equivalently, 
for the moment map to have the division property.  
We recall (slight generalizations of) these conditions in
Corollaries~\ref{l1} and \ref{l2}.  Our proof of Proposition~\ref{torus}
thus closes the gap in the proof of these Corollaries.

\begin{Proposition} \labell{torus}
Let $T$ be a torus acting linearly and symplectically 
on a symplectic vector space $V$, and let $\Phi : V \to \t^*$
be a corresponding moment map. Then for any $l \geq 0$,
the map 
\begin{equation} \labell{eq:IxPhi}
  	I \times \Phi \ : \ \R^l \times V \to \R^l \times \t^*, \quad
(u,v) \mapsto (u,\Phi(v))
\end{equation}
has the division property (see Definition \ref{property}).
\end{Proposition}

\begin{pf}
We may identify $V$ with $\C^n$ (where $n= \frac{1}{2} \dim V$) in
such a way that $T$ acts as a subtorus of the standard maximal torus
${\Bbb T}^n$ of $U(n)$.  The $T$-moment map, $\Phi$, 
is the composition of the ${\Bbb T}^n$-moment map,
$$
	F: (z_1, \ldots, z_n) \mapsto ( |z_1|^2, \ldots, |z_n|^2 ),
$$
with a linear projection $\pi : \Lie ({\Bbb T}^n)^* = \R^n \to \t^*$.
The diagram
$$ 
	\begin{array}{rcc}
   \R^l \times \C^n & \stackrel{I \times F}{\to} & \R^l \times \R^n \\
    & {\scriptstyle I \times \Phi} \searrow
      \phantom{\scriptstyle I \times \Phi} &
\phantom{\scriptstyle I \times \pi} \downarrow
	 {\scriptstyle I \times \pi} \\
    & & \R^l \times \t^*
	\end{array}
$$
commutes.  Consider the action of ${\Bbb T}^n$ on $\R^l \times \C^n$ which is
trivial on $\R^l$ and standard on $\C^n$.  It is easy to see that the
${\Bbb T}^n$-invariant polynomials on $\R^l \times \C^n$ are the polynomials
in $u_1,\ldots,u_l,|z_1|^2,\ldots,|z_n|^2$ where $u_1,\ldots,u_l$ are
coordinates on $\R^l$. That is, the ${\Bbb T}^n$-invariant polynomials on
$\R^l \times \C^n$ are the pullbacks by $I \times F$ of polynomials on
$\R^l \times \R^n$.  It follows by a theorem of G. Schwarz
\cite{sch:sm} that the ${\Bbb T}^n$-invariant smooth functions on $\R^l
\times \C^n$ are the pullbacks by $I \times F$ of smooth functions on
$\R^l \times \R^n$.

In particular, if $f \in \Cinf(\R^l \times \C^n)$ is constant on the
level sets of $I \times \Phi$, it is ${\Bbb T}^n$-invariant (because
it is also constant on the level sets of $I \times F$), so there
exists a smooth function, $g \in \Cinf(\R^l \times \R^n)$, with $g
\circ (I \times F) =f$.  Since the image of $I \times F$ is the closed
positive orthant $\R^l \times \R^n_+$, the function $g$ is constant on
the intersections of the level sets of $I \times \pi$ with this
orthant.  It is not clear {\em a priori} that
$g$ can be chosen to be constant on the level sets of $I \times \pi$.

We would like to show that $f$ is a formal pullback with respect to 
$\Phi$.  It is
enough to find, for every value $\alpha$ in the image of $I \times
\Phi$, a Taylor series on $\R^l \times \t^*$, centered at $\alpha$,
whose pullback is equal to the Taylor series of $f$ at all the points
in the level set $(I \times \Phi) \inv (\alpha)$.

Since $f$ can be written as a composition $f = (I \times F) \circ g$ for
$g \in \Cinf(\R^l \times \R^n)$, it is enough to find a formal power
series on $\R^l \times \t^*$, centered at $\alpha$, whose pullback to
$\R^l \times \R^n$ is equal to the Taylor series of $g$ at every point in
the level set $(\R^l \times (\R_+)^n) \cap (I \times \pi)^{-1} (\alpha)$.

We can choose new coordinates on $\R^l \times \R^n$ and on $\R^l
\times \t^*$ such that the map $I \times \pi$ becomes the projection
$(x_1, \ldots, x_k, y_1, \ldots, y_{m}) \mapsto (x_1, \ldots, x_k)$,
where $k= \dim T$ and $k+m = l+n$.
It is sufficient to show that for every point $p \in (I \times
\pi)^{-1}(\alpha) \cap (\R^l \times (\R_+)^n)$, the mixed partial
derivatives
$$
	\frac{\partial^{a+b} g }{\partial y^a \partial x^b} (p)
$$
are
\begin{enumerate}
\item
independent of $p$
\item
equal to zero whenever the multi-index $a$ is not zero.
\end{enumerate}
Condition 2 implies condition 1 because
any two such points $p_1$ and $p_2$ can be connected by a smooth path which
lies entirely inside the closed positive orthant and on which the
$y_i$-coordinates are constant.

Since in the interior of the positive orthant the function
$g$ is constant on level sets of the projection,
$\frac{\partial^{a+b} g}{\partial y^a \partial x^b} (q) = 0$ 
for every $q$ in the interior of the positive orthant provided $a \neq 0$.
By continuity, this also holds for all $q$ in the closed orthant.
This proves condition 2. 

We have shown that if a function $f$ is (locally) constant on the
level sets of $\Phi$ then it is a formal pullback with respect
to $\Phi$.
By Corollary \ref{criterion} of Theorem \ref{pullback}, 
the map $\Phi$ has the division property.
\end{pf}

\begin{Corollary} \labell{l1} 
A proper moment map for a Hamiltonian action of a torus on a
connected symplectic manifold has the division property.
\end{Corollary}

\begin{pf}
This was proved by Lerman in \cite{l:centralizer}.  The essential
ingredients are the facts that the level sets of $\Phi$ are connected
(cf.\ Proposition \ref{connected}) and that, by the
Marle-Guillemin-Sternberg local normal form, the moment map locally
looks like the map \eqref{eq:IxPhi}.
\end{pf}

\begin{Notation}
Denote by $\fg^*_\reg$ the elements of $\g^*$ whose stabilizers under
the coadjoint action of $G$ are tori.
\end{Notation}

\begin{Corollary} \labell{l2}
Let $G$ be a compact Lie group acting on a symplectic manifold $M$ 
with a proper moment map $\Phi: M \to \g^*$. 
Suppose that the image $\Phi (M)$ is contained in $\fg^*_\reg$.
Then $\Phi$ has the division property.
\end{Corollary}

\begin{pf}
This too was proved in \cite{l:centralizer}.
Like Corollary \ref{l1}, it follows from the fact that the level
sets of a proper moment map are connected and from the fact
that, by the local normal form, the moment map on a 
neighborhood of a point $x$ locally looks like the map 
\eqref{eq:IxPhi}, with the torus $T$ being the stabilizer of $\Phi(x)$.
\end{pf}

In section \ref{sec:SU2} we will use the following, somewhat stronger,
statement.

\begin{Corollary} \labell{nonabelian}
Let $\Phi : M \to \g^*$ be a proper moment map associated to a
Hamiltonian action of a compact connected Lie group $G$ on a connected
symplectic manifold $M$.  Then the restriction $\Phi |_
{\Phi\inv(\g^*_\reg)} : \Phi\inv(\g^*_\reg) \to \g^*_\reg$ has the
division property.
\end{Corollary}

\begin{pf}
The proof, {\em mutatis mutandis,} is the same as the proof of
Corollary~\ref{l2}.
\end{pf}
Note that the hypothesis that the moment map is proper can be replaced
by the hypotheses that it is semi-proper as a map into some open subset
of $\g^*$ and that its level sets are connected.

\begin{Remark}
For Lie groups of rank 1 or 2 we can prove Corollary~\ref{nonabelian}
without appeal to the theorem of Bierstone and Milman.  The key point
is:

\begin{Lemma}
If $V$ is a symplectic vector space and $\Phi:V \to \R$ is a moment 
map for a linear circle action, $\Phi$ has the division property.
\end{Lemma}

\begin{pf}
Assume, without loss of generality, that $\Phi(0) = 0$.  We can
identify $V$ with $\C^n$ in such a way that the circle action becomes
$\lambda \cdot (z_1, \ldots, z_n) = (\lambda^{m_1}z_1, \ldots,
\lambda^{m_n}z_n)$ and the moment map becomes $\Phi(z_1,\ldots,z_n) =
\sum m_i |z_i|^2$ for some integers $m_1,\ldots,m_n$ which are not
all zero.  Note that the level sets of $\Phi$ are connected.  Let $f$ be a
smooth function on $V$ that is constant on these level sets.  We need to
find a smooth function, $\varphi$, on $\R$ such that $f = \varphi \circ
\Phi$.

If $\text{image}(\Phi) = \R$, the exponents $m_i$ cannot all have the
same sign.  There exists then a unique function $\varphi$ on $\R$ such
that $f = \varphi \circ \Phi$.  Since the image of the regular points
of $\Phi$ is then all of $\R$, the function $\varphi$ is smooth.

Otherwise, $\text{image}(\Phi)$ is a ray and all the $m_i$'s have the
same sign.  The moment map $\Phi$ is a composition of the map $J(z_1,
\ldots, z_n) = (|z_1|^2, \ldots, |z_n|^2)$ from $V$ to $\R^n$ and of the
linear projection $\pi(x_1,\ldots,x_n) = \sum m_i x_i$ from $\R^n$ to
$\R$.  By a theorem of G.  Schwarz, 
there exists a smooth function $\tilde{\varphi}$ on $\R^n$
such that $f = \tilde{\varphi} \circ J$.  The diagonal line
$\{(t,\ldots,t)\}$ in $\R^n$ is transverse to the kernel of the
projection $\pi$; this follows from the fact that all the $m_i$'s have
the same sign.  Therefore there exists a linear map $s : \R \to \R^n$
whose image is the diagonal line and such that $\pi \circ s$ is the
identity. The function $\varphi = \tilde{\varphi} \circ s$ is smooth
on $\R$ and satisfies $f = \varphi \circ \pi$; the reason is that if
$x$ is in $\Phi(V)$ then $s(x)$ is in $J(V)$.
\end{pf}

A similar proof works for actions of tori of rank 2 but not for
higher dimensional tori. 
%The reason is that a wedge in $\R^2$
%has two extremal rays, whereas a convex polyhedral cone in $\R^k$,
%$k \geq 3$, may have more than $k$ extremal rays.
This is because if the image of the positive orthant via a linear 
projection from $\R^n$ to $\R^k$ has more than $k$ extremal rays,
this image is not equal to the image of the intersection
of an affine subspace with the positive orthant in $\R^n$.
\end{Remark}

\section{The centralizer of $SU(2)$-invariant functions.}
\labell{sec:SU2} 

Suppose that the group $SU(2)$ acts on a connected symplectic
manifold $(M,\omega)$ in a Hamiltonian fashion and that the
moment map, $\Phi : M \to su(2)^*$, is proper.
In this section we completely characterize the centralizer of the
$SU(2)$-invariant functions in $\Cinf(M)$.

\begin{Theorem} \labell{thm:SU2}
\begin{enumerate}
\item \labell{first}
If the zero level set $\Phi^{-1}(0)$ is empty, or if there is a point
in the zero level set which is not fixed by $SU(2)$, the moment map has 
the division property; equivalently, the centralizer of the invariant 
functions consists of the algebra of collective functions.
\item \labell{second}
If the zero level set $\Phi^{-1}(0)$ is nonempty, and if all the points
in this level set are fixed by $SU(2)$, the centralizer of the
invariant functions is strictly larger than the algebra of
$SU(2)$-collective functions.  The action of $SU(2)$ then extends to an
action of $U(2)$ with the same orbits and hence the same
algebra of invariant functions:
$$
 \Cinf(M)^{SU(2)} = \Cinf(M)^{U(2)}.
$$ 
The centralizer of the algebra of invariant functions consists 
of the $U(2)$-collective functions:
$$ 
(\Cinf(M)^{SU(2)})^c = \tilde{\Phi}^* \Cinf(u(2)^*),
$$
$\tilde{\Phi}$ being the $U(2)$-moment map.
\end{enumerate}
\end{Theorem}

\begin{pf}
Let $f$ be in the centralizer of the invariant functions. 
By Corollary~\ref{cor1}, $f$ pushes forward to a continuous function, 
$\Phi_*f$, on the image, $\Phi(M)$.
Recall that the function $f$ is collective if
and only if there exists a smooth function $\varphi$ on
$su(2)^*$ with $f = \varphi \circ \Phi$, i.e., if and only if
$\Phi_* f$ is smooth on the interior of $\Phi(M)$ and extends 
to a smooth function on all of $su(2)^*$.

Since the coadjoint action of $SU(2)$ factors through the
standard representation of $SO(3)$ on $su(2)^* \cong \R^3$, the
origin is the only point which has a nonabelian stabilizer.  By
Corollary~\ref{nonabelian} there exists a smooth function $\varphi$ on
$\R^3 \smallsetminus \{ 0 \}$ such that $\varphi \circ \Phi = f$ on
$\Phi ^{-1}(\R^3 \smallsetminus \{ 0 \})$.  The function $\varphi$ may
or may not extend to a function on $\R^3$ which is smooth at the origin.

If the zero level set $\Phi\inv(0)$ is empty then, since the moment
map is proper, the image $\Phi(M)$ avoids a whole neighborhood of $0$.
The push-forward, $\Phi_*f$, then extends to a smooth function on $\R^3$.

If the zero level set is nonempty, the image $\Phi(M)$ contains
a neighborhood of $0$ (unless $M$ is a single point). 

Assume that the zero level set contains a point $p$ 
whose stabilizer is not all of $SU(2)$.  
This stabilizer can be either zero or one dimensional.
If it is zero dimensional, the moment map $\Phi$ is a submersion 
from a neighborhood of $p$ to a neighborhood of $0$, 
hence the push-forward $\Phi_*f$ is smooth at $0$.

If the stabilizer of $p$ is one dimensional, it is either a maximal torus 
in $SU(2)$ or the normalizer of a maximal torus.
Denote the stabilizer by $H$.
The orbit of $p$ can be identified with the quotient $SU(2) / H$, 
which is either a sphere, $S^2$,  or a real projective plane, 
$\R \PP^2 = S^2 / \Z_2$.
It is no loss of generality to assume that the orbit is $S^2$.
By Theorem~\ref{locsymp}, a neighborhood of the orbit is equivariantly
diffeomorphic to a neighborhood of the zero section in the bundle 
$SU(2) \times_H (\h^0 \oplus V)$ over $SU(2) / H$.
Here $V$ is a symplectic vector space and $\h^0$ is the
annihilator of $\h$ in $su(2)^*$.  The cotangent bundle of the orbit
is the sub-bundle $SU(2) \times_H \h^0$. The moment map on
$SU(2) \times_H \h^0$ restricted to the fiber over $eH \in SU(2)/H$ is
the inclusion of $\h^0$ in $su(2)^*$.  The image of $SU(2) \times_H
\h^0$ under the moment map is $SU(2) \cdot \h^0 = su(2)^*$ (the
annihilator $\h^0$ is a plane through the origin in $su(2)^* \cong
\R^3$ and $SU(2)$ acts by rotations).  We are thus reduced to proving
the following:

\begin{Proposition} \labell{smile}
Consider the bundle
$TS^2 = \{ (x,y) \in \R^3 \times \R^3 \mid |x|^2=1, x \cdot y=0 \}$ 
and consider the map from it to $\R^3$ given by $(x,y) \mapsto y$.
Suppose that $\varphi$ is a continuous function on $\R^3$ that is smooth
on $\R^3 \smallsetminus \{0\}$ and whose pullback to $TS^2$ via the map
$(x,y) \mapsto y$ is smooth.  Then $\varphi$ is smooth on $\R^3$.
\end{Proposition}

The restriction of the bundle $TS^2$ to the equator $S^1 \subset S^2$
can be identified with the cylinder
$$S^1 \times \R \times \R$$
with coordinates $\theta \mod 2\pi$, $r$, and $u$.  In these
coordinates the embedding $S^1 \times \R \times \R \hookrightarrow
TS^2$ is given by 
$$
 (\theta, r, u) \mapsto (\sin \theta, - \cos \theta, 0, r \cos \theta, r
\sin \theta, u).
$$
On this cylinder the map to $\R^3$ is
\begin{equation} \labell{J}
  J : S^1 \times \R \times \R \to \R^3, \quad  \quad
  J(\theta,r,u) = (r\cos\theta, r\sin\theta,u).
\end{equation}
Its image, $J(S^1 \times \R \times \R)$, is all of $\R^3$.  
We are done if we can show

\begin{Proposition} \labell{cylinder}
Suppose that $\varphi$ is a continuous function on $\R^3$ 
that is smooth on $\R^3 \ssminus \{0\}$ and whose pullback,
$J^*\varphi$, is smooth on $S^1 \times \R \times \R$.
Then $\varphi$ is smooth on $\R^3$.
\end{Proposition}

We prove Proposition~\ref{cylinder} in a string of lemmas.  We
keep the notation of the Proposition.

\begin{Lemma} \labell{otoboos}
Suppose $\varphi$ is a smooth function on $\R^3 \ssminus \{0\}$
whose pullback, $J^*\varphi$, extends to a smooth function on 
$S^1 \times \R \times \R$.  Then for every multi-index $\alpha$ the function
$J^*{\partial^\alpha \varphi \over \partial y^\alpha}$
extends to a smooth function on $S^1 \times \R \times \R$.
\end{Lemma}

\begin{pf}
Let $f \in C^\infty(S^1 \times \R \times \R)$
be a smooth extension of $J^* \varphi$.
We can write
$$ f(\theta,r,u) = f(\theta,0,u) + r \psi(\theta,r,u)$$
where $\psi$ is smooth on $S^1 \times \R \times \R$.
For all $\theta$ we have $f(\theta,0,u) = f(0,0,u)$;
for $u \neq 0$ this holds because $f$ is then a pullback,
and for $u=0$ this holds by continuity.
So we can write
\begin{equation} \labell{train}
 f(\theta,r,u) = f(0,0,u) + r \psi(\theta,r,u).
\end{equation}
Since $$\displaystyle \dd{y_1} =
 J_*\left(\cos\theta \dd{r} - {1 \over r} \sin\theta \dd{\theta}\right),$$
since $f$ is an extension of $J^*\varphi$, and since \eqref{train}
holds, we have
\begin{equation} \labell{sm-ex}
   J^* {\partial \varphi \over \partial y_1} =
   \cos\theta {\partial (r\psi) \over \partial r}
   - \sin\theta {\partial \psi \over \partial \theta}
\end{equation}
on $J^{-1}(\R^3 \ssminus\{ 0\})$.  The right hand side of
\eqref{sm-ex} provides a smooth extension of $ J^* {\partial \varphi
\over \partial y_1} $ to all of $S^1 \times \R \times \R$.  By a
similar argument $J^* {\partial \varphi \over \partial y_2}$ and $J^*
{ \partial \varphi \over \partial y_3}$ also extend to functions in
$C^\infty(S^1 \times \R \times \R)$.

We have shown that all the first partials, $J^* {\partial \varphi \over
\partial y_i}$, extend to smooth functions on $S^1 \times \R \times \R$.
The lemma follows by a successive application
of this argument to the partial derivatives of $\varphi$.
\end{pf}

\begin{Lemma}\labell{C1}
Suppose that $\varphi$ is a continuous function on $\R^3$
which is twice continuously differentiable on $\R^3 \ssminus \{0\}$. 
Suppose that the first and second partial derivatives of $\varphi$
extend to continuous functions on $\R^3$. Then $\varphi$ is 
continuously differentiable at the origin.
\end{Lemma}

\begin{pf}
Let $g_i \in C^0(\R^3)$ be the continuous extension of $\partial
\varphi \over \partial y_i$.  We need to show that the partial
derivative $\partial \varphi \over \partial y_i$ exists at the origin
and is equal to $g_i(0)$.  By restricting attention to the appropriate
line in $\R^3$ the problem becomes one dimensional.  

An easy estimate shows that a function $\varphi$ in $C^0(\R) \cap
C^2(\R \ssminus \{0\})$ whose first derivative extends to a continuous
function on $\R$ and whose second derivative is bounded near $0$ 
is continuously differentiable at zero.
\end{pf}

\begin{Lemma} \labell{continuous}
Suppose that $\varphi$ is a continuous function on $\R^3 \ssminus \{0\}$
and that its pullback $J^*\varphi$ extends to a continuous function on the 
cylinder $S^1 \times \R \times \R$. Then $\varphi$ extends to a continuous 
function on $\R^3$.
\end{Lemma}

\begin{pf}
Let $f \in C^0(S^1 \times \R \times \R)$ be the continuous extension
of $J^* \varphi$.  For all $\theta$ we have
\begin{equation} \labell{const}
 f(\theta,0,u) = f(0,0,u)
\end{equation}
for all $u \neq 0$, because $f$ is then a pullback.
By continuity \eqref{const} also holds when $u=0$.
Hence $f$ is constant on the level sets of $J$
and descends to a function on $\R^3$.
This function coincides with $\varphi$ on $\R^3 \ssminus 0$.
It is continuous on $\R^3$ because the map $J$ is proper.
\end{pf}

\begin{pf*} {Proof of Proposition~\ref{cylinder}}
By Lemma \ref{otoboos}, for every multi-index $\alpha$ the pullback
$J^*{\partial^\alpha \varphi \over \partial y^\alpha}$ extends to a
smooth function on $S^1 \times \R \times \R$.  Therefore, by
Lemma~\ref{continuous}, all partial derivatives $\partial^\alpha
\varphi \over \partial y^\alpha$ extend to continuous functions on
$\R^3$.  In particular, the first and second partial derivatives of
$\varphi$ extend to continuous functions on $\R^3$. Hence, by
Lemma~\ref{C1}, $\varphi \in C^1 (\R^3)$.

This proves that if $\varphi \in  C^\infty(\R^3 \cap C^0(\R^3) 
\ssminus \{0\})$ and $J^*\varphi \in C^\infty (S^1 \times \R \times \R
)$, $\varphi \in C^1 (\R^3)$.

Now consider a first partial, $\partial \varphi \over \partial y_i$.
We know that it is in $C^\infty(\R^3 \ssminus \{0\}) \cap C^0(\R^3)$
and that $J^* ({\partial \varphi \over \partial y_i}|_{\R^3 \ssminus
\{0\}})$ extends to a smooth function on $S^1 \times \R \times \R$
(hence, by continuity, $J^* ({\partial \varphi \over \partial y_i}) \in
C^\infty (S^1 \times \R \times \R)$.)  Therefore, by the above argument,
$\partial \varphi \over \partial y_i$ is in $ C^1 (\R^3)$.

The argument proceeds by induction on the length of the
multi-index $\alpha$ in the partial derivative $\partial^\alpha
\varphi \over \partial y^\alpha$.
\end{pf*}

This completes the proof of Proposition~\ref{smile} and hence the
proof of part \ref{first} of Theorem~\ref{thm:SU2}.
Let us now prove part \ref{second} of the theorem.
As before, suppose that $\Phi : M \to su(2)^*$
is a proper moment map for a Hamiltonian action
of $SU(2)$ on a symplectic manifold $(M,\omega)$.
Let 
$$ Z = \Phi\inv(0) $$
be the zero level set of the moment map.
Every $SU(2)$-fixed point must lie on $Z$, because the moment map
is equivariant.  We assume that $Z$ is nonempty and coincides with
the set of fixed points.  

The normal bundle of $Z$ in $M$ is a symplectic vector bundle, hence 
can be given the structure of a complex Hermitian vector bundle. 
The key to the proof is that the representation of $SU(2)$ on a fiber 
of this
normal bundle must be the standard representation of $SU(2)$ on $\C^2$.

By the equivariant symplectic embedding theorem (see \cite{EMS} or
\cite{ma}) there exists an $SU(2)$-equivariant diffeomorphism from a
neighborhood of the zero section in the normal bundle to a neighborhood
of $Z$ in $M$ such that on a fiber of the normal bundle,
the pull-back of the moment map coincides with the moment map for the 
linear action of $SU(2)$ on that fiber.

Recall that for every non-negative integer $m$ there exists exactly
one irreducible representation of $SU(2)$ of complex dimension $m+1$.
The moment map for the representation of $SU(2)$ on $\C^{m+1}$, 
as a map from $\C^{m+1}$ into $\C \times \R \cong su(2)^*$, can
be easily computed to be:
$$
 (u_0, \ldots, u_m) \mapsto
 (\ol{u_0} u_1 + \ol{u_1} u_2 + \ldots + \ol{u_{m-1}} u_m \, ,\,
    {m \over 2} |u_0|^2 + {m-2 \over 2} |u_1|^2 + \dots + 
    {-m \over 2} |u_m|^2) .
$$
If $m$ is greater than one then the zero level set of this moment map
contains the vector $v = (1,0, \ldots, 0,1)$, which is not fixed
by the action of $SU(2)$.
Therefore, in the representation of $SU(2)$ on the fibers of the
normal bundle of $Z$ in $M$, all the irreducible components are
the standard representation on $\C^2$.  Moreover, this representation 
can occur only once: on $\C^2 \oplus \C^2$, the moment map is
$$
   ((u_0,u_1),(v_0,v_1)) \mapsto
   (\ol{u_0} u_1 + \ol{v_0} v_1 ,\,
   \half |u_0|^2 - \half |u_1|^2 + \half |v_0|^2 - \half |v_1|^2),
$$
and its zero level sets contains a vector, $((1,0),(0,1))$, 
which is not fixed by $SU(2)$.  

We have shown that the representation of $SU(2)$ 
on the fibers of the normal bundle of $Z$ in $M$
is the irreducible representation on $\C^2$ with the moment map
$$ \Phi (u_0 , u_1) = ( \ol{u_0} u_1 , \half |u_0|^2 - \half |u_1|^2 ).
$$
Note that $||\Phi||$ $=$ 
$(|\ol{u_0} u_1|^2 + {1 \over 4} (|u_0|^2 - |u_1|^2)^2)^{1/2}$ 
$ = \half ( |u_0|^2 + |u_1|^2 ) $ $=$ $ \half ||u||^2$.  
Since the square of the norm on a Hermitian vector bundle 
is a smooth function, the function $||\Phi||$ is smooth on a 
neighborhood of $Z$ in $M$. In particular, its Hamiltonian vector field 
is well defined.  We will show that this vector field is tangent to the
$SU(2)$-orbits, that it generates a circle action on $M$ which
commutes with the action of $SU(2)$, and that these actions fit
together to an effective action of $SU(2) \times_{Z_2} S^1 = U(2)$.

At points of $Z$ itself, $d||\Phi|| =0$, so the Hamiltonian vector field 
of $||\Phi||$ vanishes. 
Hence it is sufficient to prove the required properties on
the complement, $M \smallsetminus Z$.

Since the function $||\Phi||$ is $SU(2)$-invariant, its Hamiltonian
flow preserves $\Phi$. In particular, it preserves the cross section,
$\Phi^{-1}(\ell)$, where $\ell$ is an open ray in $su(2)^* \cong \R^3$
emanating from the origin.  On the cross section, the function
$||\Phi||$ coincides with the Hamiltonian generator of the action of
the circle subgroup which stabilizes the ray $\ell$.  Hence, along the
cross section, the Hamiltonian flow of the function $||\Phi||$
coincides with the action of this circle subgroup.  The intersection
of all these circles is the center, $Z_2$, of $SU(2)$.  Hence the
action of $SU(2) \times S^1$ descends to an effective action of $SU(2)
\times_{Z_2} S^1 \cong U(2)$.

To prove that the centralizer of the invariants consists of the 
$U(2)$-collective functions we need, by Theorem \ref{duality},
to prove that $\Cinf(M)^\Phi = \tPhi^* \Cinf(u(2)^*)$.
Note that the image under $\Phi$ of the normal bundle of $Z$
in $M$ coincides with the image of one fiber, and that
the level sets of the $SU(2)$-moment map, $\Phi$, coincide with the
level sets of the $U(2)$-moment map, $\tilde{\Phi}$.  Hence it is
sufficient to show that the $U(2)$ moment map on $\C^2$,
$$
F : \C^2 \to \C \times \R^2, \quad
  (u_0,u_1) \mapsto (\ol{u_0} u_1 , |u_0|^2 , |u_1|^2),
$$
has the division property, i.e., that it satisfies
\begin{equation} \labell{end}
  \Cinf(\C^2)^{F} = F^* \Cinf(\C \times \R^2).
\end{equation}
The functions in $\Cinf(\C^2)^{F}$ are exactly
the $S^1$-invariant functions on $\C^2$. The coordinates of
$F$ generate the ring of $S^1$-invariant polynomials
on $\C^2$. Schwarz's theorem \cite{sch:sm} then gives
\eqref{end}.
This completes the proof of Theorem \ref{thm:SU2}.
\end{pf}

\appendix
\section{Local normal form for the moment map and implications}
\labell{sec:normal}

In this section we describe a version of the Marle-Guillemin-Sternberg
normal form of the moment map (Theorems \ref{cross-section} and
\ref{locsymp}, Corollary~\ref{use_cross-section}) and some of its
implications (Lemmas~\ref{smooth-curve} and~\ref{local-properties},
Remark~\ref{cone}) which were use in the proofs of
Theorems~\ref{duality} and~\ref{pullback}.  This version of the
Marle-Guillemin-Sternberg normal form describes the moment map for a
Hamiltonian action of a compact Lie group $G$ on a symplectic manifold
$(M,\omega)$ in terms of nice ``coordinates'' on a neighborhood of an
orbit, $G \cdot x$, in $M$.  The construction of this normal form is
carried out in two steps: the construction of the symplectic cross-section
(Theorem \ref{cross-section} and Corollary~\ref{use_cross-section}),
and the construction of the normal form near an isotropic orbit (Theorem
\ref{locsymp}).

\begin{Theorem}[The symplectic cross-section]
	\labell{cross-section}
Let $\Phi : M \to \g^*$ be a moment map associated to a Hamiltonian
action of a compact Lie group $G$ on a symplectic manifold $(M,\omega)$.
Let $\alpha$ be a point in $\g^*$ and let $G_\alpha$ denote its stabilizer
under the coadjoint action. 

Then for a sufficiently small slice  $S$ at $\alpha$ for the coadjoint
action of $G$,  the preimage $R:= \Phi \inv (S)$ is a {\em symplectic}
$G_\alpha$-invariant submanifold of $M$.  Moreover, the action of
$G_\alpha$ on $R$ is Hamiltonian and the restriction of the moment map
$\Phi$ to $R$ followed by the natural projection $\fg^* \to
\fg^*_\alpha$ is a corresponding  moment map.
\end{Theorem}
\begin{pf}
Theorem~26.7 in \cite{g-s:sympl-tech}.
\end{pf}
\begin{Remark}
The submanifold $R$ is called a {\em symplectic cross-section}.
\end{Remark}

Since $S$ is a slice, the open invariant neighborhood $G\cdot S\subset
\fg^*$ of the coadjoint orbit of $\alpha$ 
is equivariantly diffeomorphic to the associated
bundle $G\times _{G_\alpha} S$.  Since the moment map $\Phi$ is
equivariant, $G\cdot R$ is an open invariant subset of the manifold
$M$ and it is equivariantly diffeomorphic to the associated bundle
$G\times _{G_\alpha} R$.  Up to these identifications, the
moment map $\Phi : G\cdot R\to G\cdot S$ is the equivariant bundle map
$G\times _{G_\alpha} R \to G\times _{G_\alpha} S,$  $[g, r] \mapsto
[g, \Phi_R (r)]$, where $\Phi _R$ is the restriction of $\Phi$ to the
cross-section $R$, and $[g,r]$ and $[g, \Phi_R(r)]$ denote the
$G_\alpha$-orbits of $(g,r) \in G\times R$ and $(g,\Phi_R(r))\in
G\times S$   respectively.  If the slice $S$ is sufficiently small or,
more generally, is well-chosen, then the restriction to $S$ of the
projection $\fg^* \to \fg^*_\alpha$ is a diffeomorphism onto the
image.  This proves the following:

\begin{Corollary} \label{use_cross-section} 
Let $\Phi :M\to \fg^*$ be as in Theorem~\ref{cross-section},  $x$
be a point in $M$, $\alpha = \Phi (x)$, $S$ a small slice at $\alpha$
and $R = \Phi\inv (S)$ the corresponding cross-section.  Then, up to
composition with diffeomorphisms, the moment map $\Phi$ in an
invariant neighborhood of $x$ is a bundle map
$$
G\times _{G_\alpha} R \to G\times _{G_\alpha} S, \quad\quad 
[g, r] \mapsto [g, \Phi_R (r)],
$$
where $\Phi_R$ is the $G_\alpha$ moment map (up to a composition with
a diffeomorphism).
\end{Corollary}

\begin{Remark} \label{iso_remark}
Since the point $\alpha$ is fixed by the action of $G_\alpha$, 
the $G_\alpha$-orbit through any point $x\in \Phi\inv
(\alpha)$ is {\em isotropic} in the symplectic cross-section, $R$.
Since the orbit $G_\alpha \cdot x$ is isotropic, the tangent space $T
_x(G_\alpha \cdot x)$ is contained in its symplectic perpendicular $T
_x(G_\alpha \cdot x)^{\omega_R}$.  Hence the quotient $V = T _x(G_\alpha
\cdot x)^{\omega_R}/ T_x(G_\alpha \cdot x)$ is a symplectic vector space.
It is called the {\em symplectic slice \/} at $x$.  Note that $V$
is a natural symplectic representation of the isotropy group $G_x$ of $x$.
\end{Remark}

\begin{Theorem}[Local normal form near an isotropic orbit]
\labell{locsymp} 
Let $\Psi : N \to \fh^*$ be a moment map associated to a Hamiltonian
action of a compact Lie group $H$ on a symplectic manifold
$(N,\sigma)$.  Suppose the orbit $H\cdot x$ is {\em isotropic} in $N$.
Let $H_x$ denote the stabilizer of $x$ in $H$,
let $\fh_x^0$ denote the annihilator of its Lie algebra in $\fh^*$,
and let $H_x \to Sp (V)$ denote the symplectic slice representation.

Given an $H_x$-equivariant embedding, $i:\fh_x^* \to \fh^*$,  
there exists an $H$-invariant closed two-form, $\omega_Y$, on the manifold
$Y= H\times_{H_x} \left(\fh_x^\circ \times V \right)$, such that
\begin{enumerate}
\item the form $\omega_Y$ is nondegenerate near the zero section of the
bundle $Y\to H/H_x$, 
\item 
a neighborhood $U_x$ of the orbit of $x$ in $N$ 
is equivariantly symplectomorphic
to a neighborhood of the zero section in $Y$, and
\item 
the action of $H$ on $(Y, \omega _Y)$ is Hamiltonian with a moment map
$\Phi_Y : Y \to \fh^*$ given by
$$ 
\Phi _Y ([g, \eta, v])= Ad^\dagger (g)\left(\eta + i( \Phi_{V} (v))\right)
$$ 
where $Ad^\dagger$ is the coadjoint action, and $\Phi_{V} : V\to \fh_x^*$
is the moment map for the slice representation of $H_x$.
\end{enumerate}
Consequently, the equivariant embedding $\iota: U_x \hookrightarrow Y$
intertwines the two moment maps, up to translation:
$\Psi |_{U_x} = \Phi_Y \circ \iota + \Phi(x)$.
\end{Theorem}

\begin{pf}
Theorem~\ref{locsymp} is essentially an equivariant version of Weinstein's
isotropic embedding theorem. See \cite{g-s:normal} or \cite{g-s:sympl-tech}.
\end{pf}

One easy consequence of the local normal form is the following
technical lemma, which was used in the proof of Theorem~\ref{duality}.

\begin{Lemma}\labell{smooth-curve}
Let $\Phi : M \to \g^*$ be a moment map associated to a Hamiltonian action
of a compact Lie group, $G$, on a symplectic manifold, $(M,\omega)$.
Then any two points in a connected component of a level set of $\Phi$
can be joined by a piece-wise smooth curve that lies in the level set.
\end{Lemma}

\begin{pf} It is sufficient to prove the theorem for a neighborhood of a 
point $x$ in $M$.  By Theorem~\ref{cross-section},
Corollary~\ref{use_cross-section}, and Remark~\ref{iso_remark} we may
assume that the orbit $G\cdot x$ is isotropic. By
Theorem~\ref{locsymp} it is sufficient to prove the lemma for the 
zero level set of the map $\Phi _Y$ of the local model.  
This level set is 
$G \times_{G_x} (\{0\} \times \Phi_V^{-1}(0))$.  The moment map on the
symplectic vector space $V$ is homogeneous, therefore $\Phi_V^{-1}(0)$
is a cone in $V$.  Since any point in the cone can be connected to the
vertex by a straight line, any two points in the level set can be
connected by a piece-wise smooth curve.
\end{pf}

A non-trivial consequence of the local normal form theorem are the
following lemma and remark, which were used in the proof of
Theorem~\ref{pullback}.  They essentially say that the image under the
moment map of a small invariant neighborhood of an orbit $G \cdot x$ does not
change as $x$ varies along a connected component of the level set
$\Phi\inv(\Phi(x))$:

\begin{Lemma} \labell{local-properties}
Let $\Phi : M \to \g^*$ be a moment map associated to an action of a
compact connected Lie group $G$ on a connected symplectic manifold
$(M,\omega)$.  Let $x$ be a point in $M$. 

Let $T$ be a maximal torus of the isotropy group of $\Phi (x)$ (hence
of $G$), $\t$ its Lie algebra, identified with its dual $\t^*$ and
embedded in $\g^*$ via a choice of an ${\text Ad}$-invariant inner
product, and let $\t^*_+ \subset \g^*$ be the Weyl chamber containing
$\Phi(x)$.

There exists a rational polyhedral cone $C_x$ in $\t^*$ with vertex at
$\Phi(x)$, such that for every $G$-invariant neighborhood $U$ of $x$
in $M$ there exists an $\text{Ad}^\dagger(G)$-invariant neighborhood
$\calV$ of $\Phi(x)$ in $\g^*$ such that $\Phi(U) \cap \t_+^*$
contains $C_x \cap \calV$.  Hence $\Phi (U \cap \Phi \inv (\calV))
\cap \t_+^* = C_x \cap \calV$. Moreover, if
the neighborhood $U$ is sufficiently small then for any point $y$ in
$U \cap \Phi\inv(\Phi(x))$, the cone $C_y$ is equal to the cone $C_x$.
\end{Lemma}

If the group $G$ is abelian, the result is easy.  The general case is
due to Sjamaar \cite{Sj}.  Lerman, Meinrenken, Tolman, and Woodward
later found a more ``elementary'' proof: the Lemma follows from the
proofs of Theorems 6.1 and 6.2 of \cite{LMTW}.

\begin{Remark} \labell{cone}
If $\Phi(x) = 0$, the cone $C_x$ is $\Phi_Y(Y) \cap \t^*_+$
where $Y$ is the local model described in Theorem~\ref{locsymp}. 
\end{Remark}

\end{document}